\documentclass[a4paper]{paper}

\usepackage[utf8x]{inputenc}
\usepackage{mathtools}
\usepackage{etoolbox}
\usepackage{amsfonts,mathrsfs,amssymb}
\usepackage{xspace}
\usepackage{graphicx}
\usepackage[labelformat=simple]{caption}
\usepackage{hyperref}
\hypersetup{%
    bookmarks=true,         
    pdfmenubar=true,       
    pdffitwindow=false,     
    pdfstartview={FitH},    
    colorlinks=true,       
    linkcolor=red,          
    citecolor=red,        
    filecolor=magenta,      
    urlcolor=cyan,           
    pdfborder = {0,0,0}
}

\usepackage{cleveref}
\usepackage{todonotes}
\usepackage{nicefrac}
\usepackage{paralist}

\usepackage{mymathdefs}
\usepackage{mytheorems}

\newcommand*{\FWDOP}{{\OP{T}}}

\newcommand*{\DATAFUNC}{{\OP{D}}}

\newcommand*{\MODULAR}{\ensuremath{\rho_p}}

\newcommand*{\MOREAU}{\ensuremath{M}}

\newcommand*{\TVMINL}{\ensuremath{f^{(\mathrm{l})}}}

\DeclareMathOperator{\TGV}{TGV}

\crefname{equation}{}{}
\Crefname{equation}{}{}
\crefname{item}{}{}
\Crefname{item}{}{}
\crefname{mytheorem}{theorem}{theorems}
\Crefname{mytheorem}{Theorem}{Theorems}
\crefname{mylemma}{lemma}{lemmas}
\Crefname{mylemma}{Lemma}{Lemmas}

\title{Total variation regularization with variable Lebesgue prior}
\author{Holger Kohr}

\begin{document}
\maketitle

\begin{abstract}
  This work proposes the variable exponent Lebesgue modular as a replacement for the 1-norm in total variation (TV) regularization.
  It allows the exponent to vary with spatial location and thus enables users to locally select whether to preserve edges or smooth intensity variations.
  In contrast to earlier work using TV-like methods with variable exponents, the exponent function is here computed offline as a fixed parameter of the final optimization problem, resulting in a convex goal functional.
  The obtained formulas for the convex conjugate and the proximal operators are simple in structure and can be evaluated very efficiently, an important property for practical usability.
  Numerical results with variable $L^p$ TV prior in denoising and tomography problems on synthetic data compare favorably to total generalized variation (TGV) and TV.
\end{abstract}

\section{Introduction}
\label{sec:intro}

Total variation (TV) regularization is nowadays an established technique in imaging.
It consists in finding
\begin{equation*}
  f_\lambda = \ARGMIN_{f \in \SPCX} \left[ \DATAFUNC(f) + \lambda \TV(f) \right]
\end{equation*}
for a data discrepancy measure $\DATAFUNC$, a regularization parameter $\lambda > 0$ and the TV functional
\begin{equation*}
  \TV(f) = \sup \SET[\Big]{\int_\Omega f \DIV\phi\, \D x \GIVEN \phi \in L^\infty(\Omega, \RR^d),\ \NORM{\phi}_\infty < 1}
\end{equation*}
defined on the space $\BV(\Omega)$ of functions of bounded variation on $\Omega \subset \RR^d$.
This model has proved to possess favorable properties when edge preservation and effective noise suppression are primary goals in image reconstruction.
There exist very efficient implementations of numerical solvers for TV-regularized inverse problems, which also makes the technique accessible for daily practical use in engineering and medical imaging.
Despite its success, the classical TV model has also exhibited its limitations over time, most notably the staircasing effect \cite{ring_structural_2000}, i.e. the preference of the regularizer for piecewise constant functions.

Several extensions of total variation have been proposed to tackle the problem of staircasing, with the latest contribution being total generalized variation (TGV) \cite{bredies_total_2010}, a higher-order TV variant.
It has shown to perform significantly better than TV on images with both edges and gradual intensity changes, and it will serve as a benchmark for the method proposed in this paper.
On the downside, TGV requires handling of higher order tensors leading to higher numerical complexity, and it introduces an additional regularization parameter, which makes their selection more complex.

This paper takes a different approach in varying classical TV regularization based on the observation that functionals of the form $\NORM{\nabla \cdot}_p^p$ lead to very different behavior of minimizers of $f \mapsto \DATAFUNC(f) + \lambda \NORM{\nabla f}_p^p$ for $p = 1$ and $p > 1$, with a sharp transition between the two cases:
While the case $p = 1$ corresponds to TV regularization, any exponent $p > 1$ leads to a smooth solution \cite{burger_guide_2013}.
Therefore, a natural approach would be to select the exponent $p$ \emph{dependent on spatial location}, taking the value $p = 1$ where edges should be enhanced, and $p > 1$ otherwise.

The idea of using a spatially varying exponent in a TV-like regularization method has been proposed as early as 1997 \cite{blomgren_total_1997} and put into practice in 2006 \cite{chen_variable_2006}.
Both papers as well as subsequent articles try to tackle variants of the problem
\begin{equation*}
  f^* = \ARGMIN \left[ \DATAFUNC(f) + \lambda \int_\Omega \ABS{\nabla f(x)}^{p(\nabla f(x))}\, \D x \right],
\end{equation*}
where the exponent depends directly on the image, e.g.,
\begin{equation*}
  p(\nabla f) = 1 + \frac{1}{1 + k \ABS{G_\sigma \ast \nabla f}^2}.
\end{equation*}
For the resulting non-convex problem, only partial solutions have been proposed, e.g., for a smoothed version of the integrand with a weak notion of solution \cite{wunderli_time_2010}.

A parallel track that has recently received growing interest are variable exponent \emph{sequence spaces} and their application in imaging problems.
The sequence usually represents coefficients in wavelet-type decompositions, and a sequence of exponents serves to define Besov-type spaces in this case \cite{jia_bayesian_2016}.
A very recent paper studies convergence and stability properties of a variable exponent prior in this context \cite{lorenz_flexible_2017}.
That approach has no direct connection to the variable exponent functionals considered here, although the variable Lebesgue spaces can be characterized in terms of wavelet coefficients \cite{samashvili_wavelet_2015}.

The method proposed in this paper computes the exponent $p$ in an offline step and keeps it as a fixed parameter in the then convex minimization problem.
There are two natural types of imaging problems where this approach can be applied:
\begin{inparaenum}[(1)]
  \item single-channel imaging where first the exponent is computed from the given data and then applied as prior in the subsequent minimization problem (``bootstrapping'') and
  \item dual-channel imaging where the secondary channel provides the exponent map that is used for regularization of the primary channel.
\end{inparaenum}

We start by introducing the variable Lebesgue modular, the central functional in this work.
The following section studies the properties of the modular in the light of convex optimization and derives expressions for the building blocks needed to solve the TV minimization problem introduced thereafter.
A short section presents the method for computing the exponent map, before numerical results are shown.

\section{The variable Lebesgue modular}
\label{sec:modular}

Let $\Omega \subset \RR^d$ be an open and connected domain.
Let further $p \colon \Omega \to \CINTERV{1}{\infty}$ be a measurable exponent function.
The variable Lebesgue modular \cite{cruz-uribe_variable_2013} or \emph{$p$-modular} of a function $f \colon \Omega \to \RR$ is defined as
\begin{equation}
  \MODULAR(f) \DEFEQ \int_{\Omega^*} \ABS{f(x)}^{p(x)}\, \D x + \NORM{f_\infty}_\infty,
  \label{eq:varlp:def_modular_general}
\end{equation}
where $\Omega^* = \SET{x \in \Omega \GIVEN p(x) < \infty}$ and $f_\infty$ is the restriction to $\Omega_\infty = \Omega \setminus \Omega^*$.
For vector-valued functions $F\colon \Omega \to \RR^m$, we use the same notation, in which case $\ABS{F(x)}$ stands for the Euclidean norm in $\RR^m$.
Occasionally we denote this Euclidean norm with $\ABS{}_2$ if we deem it necessary for notational clarity.

We are mainly interested in exponents between 1 and 2, and therefore we assume that $p \colon \Omega \to \CINTERV{1}{2}$ and thus $\Omega = \Omega^*$.
As a direct consequence, the modular simplifies to
\begin{equation}
  \MODULAR(f) = \int_\Omega \ABS{f(x)}^{p(x)}\, \D x.
  \label{eq:varlp:def_modular}
\end{equation}

\begin{myremark}
  Although the $p$-modular is a direct generalization of the $p$-th power of the norm $\NORM{}_p$ for constant exponent, it is not homogeneous and thus cannot be easily used to define a norm.
  To achieve this, the construction
  \begin{equation}
    \NORM{f}_p \DEFEQ \inf\SET{\lambda > 0 \GIVEN \MODULAR(f / \lambda) < 1}
    \label{eq:varlp:def_norm}
  \end{equation}
  can be used.
  The corresponding space $L^p(\Omega)$ of functions where this norm is finite, as well as related Sobolev spaces, have been studied thoroughly \cite{diening_lebesgue_2011,cruz-uribe_variable_2013}.
  Recently also spaces of functions with bounded variation have attracted attention in the context of PDEs \cite{harjulehto_minimizers_2008} and imaging \cite{harjulehto_critical_2013}.
\end{myremark}

\section{Convex optimization with variable Lebesgue modular}
\label{sec:optim}

We now consider the $p$-modular as a functional on $L^2(\Omega)$, where we assign the value $\infty$ to functions in $L^2(\Omega) \setminus L^p(\Omega)$.
This identification is justified by the fact that for bounded exponent functions, the set of finite modular and finite norm are identical \cite[Proposition 2.12]{cruz-uribe_variable_2013}.

In accordance with optimization literature, we denote by $\Gamma_0(\SPCX)$ the set of \emph{proper}, \emph{convex} and \emph{lower semicontinuous} functionals on $\SPCX$.
These minimum requirements are satisfied by the $p$-modular as a consequence of the integrand possessing these properties.

\begin{myproposition}
  The $p$-modular is proper, convex and lower semicontinuous as a functional $\MODULAR \colon \SPCX \to \CINTERV{0}{\infty}$.
\end{myproposition}

\noindent
It is similarly straightforward to see why this functional suffers from the same non-differentiability issue at 0 as the 1-norm, unless the exponent is equal to 2 everywhere.

\begin{myproposition}
  The gradient of the $p$-modular is given by
  \begin{equation*}
    \nabla \rho_p(f) = p\, \ABS{f}^{p-2}\, f.
  \end{equation*}
  It is well-defined for all functions $f$ that are nonzero at almost all points $x$ where $p(x) \neq 2$.
\end{myproposition}

As a consequence, smooth optimization cannot be applied unless the functional itself is replaced by a smoothed version, e.g., $\rho_p(\ABS{f}_\epsilon)$, using the smoothed norm $\ABS{f}_\epsilon = (\ABS{f}^2 + \epsilon^2)^{1/2}$.

We do not follow this route and instead resort to non-smooth convex optimization.
In the following we define the basice notions of convex conjugate, Moreau envelope and proximal operator.
Note that all definitions can be formulated for general Banach spaces, but we focus on the special case of the Hilbert space $L^2$.

\begin{mydefinition}
  Let $\SPCX = L^2(\Omega, \RR^m)$ with $m = 1$ for the scalar case and $m > 1$ for the vectorial case.
  Let further $\OPS \in \Gamma_0(\SPCX)$.
  The \emph{convex conjugate} $\OPS^* \colon \SPCX \to \LINTERV{-\infty}{\infty}$ of $\OPS$ is the functional
  \begin{equation}
    \OPS^*(f) \DEFEQ \sup_{g \in \SPCX} \big[ \INNER{f}{g}_\SPCX - \OPS(g) \big],\quad f \in \SPCX.
    \label{eq:optim:convex_conj}
  \end{equation}
  The \emph{Moreau envelope} or \emph{Moreau-Yosida regularization} $\MOREAU_\tau (\OPS) \colon \SPCX \to \LINTERV{-\infty}{\infty}$ of $\OPS$ is defined as
  \begin{equation}
    \MOREAU_\tau (\OPS)(f) \DEFEQ \inf_{g \in \SPCX} \left[ \OPS(g) + \frac{1}{2\tau} \NORM{f - g}_2^2 \right].
  \end{equation}
  The \emph{proximal operator} $\PROX_{\tau \OPS} \colon \SPCX \to \SPCX$ of $\OPS$ is the mapping
  \begin{equation}
    \PROX_{\tau \OPS}(f) \DEFEQ \ARGINF_{g \in \SPCX} \left[ \OPS(g) + \frac{1}{2\tau} \NORM{f - g}_2^2 \right].
  \end{equation}
\end{mydefinition}

We derive these three fundamental mappings for the functional $\MODULAR$ in the following.
As an important tool, we state a simplified version of a theorem that allows to rephrase the optimization of an integrating functional as the pointwise optimization of its integrand.

\begin{mytheorem}{\cite[Theorem 14.60]{rockafellar_variational_1998}}
  Let $\SPCX = L^2(\Omega, \RR^m)$ and $F \colon \Omega \times \RR^m$ be a normal integrand, i.e., a function whose epigraphical mapping $E(x) = \mathrm{epi}F(x, \cdot)$ is closed-valued and measurable.
  Then, the integral functional
  \begin{equation*}
    \OPI_F \colon \SPCX \to \LINTERV{-\infty}{\infty}, \quad \OPI_F(u) \DEFEQ \int_\Omega F\big(x, u(x)\big)\,\D x
  \end{equation*}
  has the property that minimization and integration can be exchanged, i.e.,
  \begin{equation}
    \label{eq:optim:pointwise_inf}
    \inf_{u\in\SPCX} \OPI_F(u) = \int_\Omega \left[\inf_{z\in\RR^m} F(x, z) \right]\, \D x,
  \end{equation}
  as long as $\OPI_F \not\equiv +\infty$ on $\SPCX$.
  Likewise, unless this common value is $-\infty$, the minimizers of the functional satisfy
  \begin{equation}
    \label{eq:optim:pointwise_arginf}
    \bar u \in \ARGINF_{u\in\SPCX} \OPI_F(u) \quad\Longleftrightarrow\quad \bar u(x) \in \ARGINF_{z\in\SPCX} F(x, z)\quad \text{for almost all } x \in \Omega.
  \end{equation}
  \label{thm:optim:pointwise}
\end{mytheorem}

The condition on the integrands to be normal is rather straightforward to check in the cases considered below.
According to \cite[Example 14.29]{rockafellar_variational_1998}, all Carath\'{e}ordory integrands, i.e., functions that are measurable in the first variable and continuous in the second, are normal integrands.
This characterization applies to all integrands in this paper.

\begin{myremark}
  Some of the subsequent results can be pieced together from the literature, e.g., the convex conjugate of the norm, the Huber function \cite{bauschke_convex_2011} or expressions for the proximal operator \cite{combettes_proximal_2011} in certain cases.
  However, for the sake of completeness, and due to the coherent kind of argumentation, we provide complete proofs of all cases.
\end{myremark}

\begin{mylemma}
  The convex conjugate of $\MODULAR$ is given by
  \begin{equation}
    \MODULAR^*(f) = \int_\Omega R\big( f(x), p(x) \big)\, \D x,
    \label{eq:optim:modular_cconj}
  \end{equation}
  where
  \begin{equation}
    R(z, p) \DEFEQ
    \begin{cases}
      \iota_{\ABS{} \leq 1}(z),                                             & \text{if } p = 1, \\[1ex]
      \ABS{z}^2 / 4,                                                        & \text{if } p = 2, \\[1ex]
      \ABS{z}^{\frac{p}{p-1}} \big[ p^{-1/(p-1)}  - p^{-p/(p-1)} \big],     & \text{else},
    \end{cases}
    \label{eq:optim:modular_cconj_integrand}
  \end{equation}
  for $z \in \RR^m$ and $p \in \CINTERV{1}{2}$.
  Here, $\iota_{\ABS{} \leq 1}$ denotes the \emph{indicator function} of the unit ball in $\RR^m$, i.e., the function that takes the value 0 in the set and $\infty$ outside.
  \label{lemma:optim:modular_cconj}
\end{mylemma}

\begin{myproof}
  Since $\sup A = -\inf (-A)$ for any set $A \subset \RR$, \Cref{thm:optim:pointwise} holds also for $\sup$ instead of $\inf$.
  For the $p$-modular, the functional to be maximized can be written as
  \begin{equation*}
    \INNER{f}{g}_\SPCX - \MODULAR(g) = \int_\Omega \left[ \INNER{f(x)}{g(x)}_{\RR^m} - \ABS{g(x)}^{p(x)} \right]\D x,
  \end{equation*}
  which means that \Cref{eq:optim:pointwise_inf} can be applied, i.e., integration and optimization can be exchanged:
  \begin{equation*}
    \MODULAR^*(f) = \sup_{g \in \SPCX} \big[ \INNER{f}{g}_\SPCX - \MODULAR(g)\big] = \int_\Omega \sup_{y \in \RR^m} \left[ \INNER{f(x)}{y}_{\RR^m} - \ABS{y}^{p(x)} \right]\D x.
  \end{equation*}
  Hence, the optimization of the integrand can be performed for each point $x \in \Omega$ separately.
  In what follows, we consider $p = p(x)$ to be a constant and set $z = f(x) \in \RR^m$.
  The new problem is to find
  \begin{equation*}
    r^*(z) = \sup_{y \in \RR^m} \big[ \INNER{z}{y}_{\RR^m} - \ABS{y}^{p} \big],
  \end{equation*}
  which is nothing but the convex conjugate of $r = \ABS{}^p$.
  Apparently, $r^*(0) = 0$.
  For other points, using polar coordinates, this problem can be rephrased as
  \begin{equation*}
    r^*(z) = \sup_{\alpha \geq 0}\sup_{\theta \in S^{m-1}} \big[ \alpha \INNER{z}{\theta}_{\RR^m} - \alpha^{p} \big] = \sup_{\alpha \geq 0} \big[ \alpha(\ABS{z} - \alpha^{p-1}) \big].
  \end{equation*}
  If $p=1$, the expression in brackets is only bounded from above if $\ABS{z} \leq 1$, in which case the supremum is $0$.
  Thus, the convex conjugate is $r^* = \iota_{\ABS{} \leq 1}$, the indicator function of the unit ball.

  Otherwise, we can differentiate $v(\alpha) = \alpha\ABS{z} - \alpha^p$ with respect to $\alpha$ and get
  \begin{equation*}
    v'(\alpha) = \ABS{z} - p\alpha^{p-1}, \quad v''(\alpha) = -p(p-1)\alpha^{p-2},
  \end{equation*}
  which shows that $\bar\alpha = (\ABS{z} / p)^{1/(p-1)}$ is a critical value of $v$, and further that $v$ is strictly concave and the critical value is the only local maximum.
  The maximum value is
  \begin{equation*}
    v(\bar \alpha) = \ABS{z}\left( \frac{\ABS{z}}{p} \right)^{\frac{1}{p-1}} - \left( \frac{\ABS{z}}{p} \right)^{\frac{p}{p-1}} = \ABS{z}^{\frac{p}{p-1}} \left[ p^{-1/(p-1)}  - p^{-p/(p-1)} \right],
  \end{equation*}
  which simplifies to $v(\bar \alpha) = \ABS{z}^2 / 4$ for $p = 2$.
  Since this value is larger than $v(0) = 0$, it is also the global maximum.
\end{myproof}

For the Moreau envelope, we can largely use the same technique to find pointwise expressions.

\begin{mylemma}
  The Moreau envelope of $\MODULAR$ can be written as
  \begin{equation}
    \MOREAU_\tau(\MODULAR)(f) = \int_\Omega T_\tau\big( f(x), p(x) \big)\, \D x,
    \label{eq:optim:modular_moreau}
  \end{equation}
  where
  \begin{equation}
    T_\tau(z, p) \DEFEQ
    \begin{cases}
      \frac{\ABS{z}^2}{2\tau},                                                      & \text{if } p = 1 \text{ and } \ABS{z} \leq \tau,  \\[1ex]
      \ABS{z} - \frac{\tau}{2},                                                     & \text{if } p = 1 \text{ and } \ABS{z} \geq \tau,  \\[1ex]
      \frac{\ABS{z}^2}{1 + 2\tau},                                                  & \text{if } p = 2,                                 \\[1ex]
      (\ABS{z} - \bar\alpha) \frac{2\bar\alpha + p(\ABS{z} - \bar\alpha)}{2\tau p}, & \text{ else}.
    \end{cases}
    \label{eq:optim:modular_moreau_integrand}
  \end{equation}
  Here, $\bar\alpha = \bar\alpha(\ABS{z}, p, \tau)$ denotes the unique solution $0 \leq \bar\alpha < \ABS{z}$ to the equation
  \begin{equation}
    \alpha + \tau p \alpha^{p-1} = \ABS{z}.
    \label{eq:optim:moreau_implicit}
  \end{equation}
  \label{lemma:optim:modular_moreau}
\end{mylemma}

\begin{myproof}
  We have, again due to \Cref{eq:optim:pointwise_inf}, that
  \begin{align*}
    \MOREAU_\tau(\MODULAR)(f)
    &= \inf_{g \in \SPCX} \int_\Omega \left[ \ABS{g(x)}^{p(x)} + \frac{1}{2\tau} \ABS{f(x) - g(x)}^2 \right] \, \D x \\
    &= \int_\Omega \inf_{y \in \RR^m} \left[ \ABS{y}^{p(x)} + \frac{1}{2\tau} \ABS{f(x) - y}^2 \right] \, \D x.
  \end{align*}
  Using the same notation as in the proof of \Cref{lemma:optim:modular_cconj}, we can rewrite the inner optimization problem as the the determination of
  \begin{align}
    \MOREAU_\tau\big(\ABS{}^p\big)(z)
    &= \inf_{y \in \RR^m} \left[ \ABS{y}^{p} + \frac{1}{2\tau} \ABS{z - y}^2 \right] \notag\\
    &= \inf_{\alpha \geq 0}\inf_{\theta \in S^{m-1}} \left[ \alpha^{p} + \frac{1}{2\tau} \big(\ABS{z}^2 - 2 \alpha\INNER{z}{\theta}_{\RR^m} + \alpha^2\big) \right]
    \label{eq:optim:moreau_sphere_min} \\
    &= \inf_{\alpha \geq 0} \left[ \alpha^{p} + \frac{1}{2\tau} \big(\ABS{z}^2 - 2 \alpha\ABS{z} + \alpha^2\big) \right] \notag\\
    &= \inf_{\alpha \geq 0} \underbrace{\left[ \alpha^{p} + \frac{1}{2\tau} \big(\ABS{z} - \alpha\big)^2 \right]}_{\EQDEF v(\alpha)}. \notag
  \end{align}
  Differentiation of $v$ yields
  \begin{align*}
    v'(\alpha) &= p \alpha^{p-1} + \tau^{-1}(\alpha - \ABS{z}), \\
    v''(\alpha) &=
    \begin{cases}
     \tau^{-1},                             &   \text{if } p = 1,   \\
     p (p - 1) \alpha^{p-2} + \tau^{-1},    &   \text{else}.
    \end{cases}
  \end{align*}
  Thus the only critical point of $v$ is given as the solution to the implicit equation
  \begin{equation*}
    \alpha + \tau p \alpha^{p-1} = \ABS{z},
  \end{equation*}
  which is the same as \Cref{eq:optim:moreau_implicit}.
  In the special case $p = 2$, this equation is easily solved for $\alpha$, resulting in $\bar\alpha = \ABS{z} / (1 + 2\tau)$.
  The positivity of $v''$ lets us conclude that this point is the only local minimum of $v$.
  For the minimum value, we get
  \begin{equation*}
    v(\bar \alpha) = \frac{\ABS{z}^2}{(1 + 2\tau)^2} + \frac{1}{2\tau}\left(\ABS{z} -  \frac{\ABS{z}}{1 + 2\tau}\right)^2 = \frac{\ABS{z}^2 + \ABS{z}^2(2\tau)}{(1 + 2\tau)^2} = \frac{\ABS{z}^2}{1 + 2\tau}.
  \end{equation*}
  This value is clearly smaller than $v(0) = \ABS{z}^2 / (2\tau)$, therefore $\bar\alpha$ is also the global minimum.

  For $p = 1$, on the other hand, the critical point lies at $\bar\alpha = \ABS{z} - \tau$.
  Here we must distinguish two cases:
  \begin{enumerate}
    \item[$\ABS{z} < \tau$:]
    The critical point is negative and thus not eligible as minimum.
    To find the global minimum, we observe that
    \begin{equation*}
      v'(\alpha) = \underbrace{1 - \tau^{-1}\ABS{z}}_{> 0} + \tau^{-1}\alpha > 0,
    \end{equation*}
    from which immediately follows that the global minimum is at $\alpha = 0$, with minimum value $v(\bar \alpha) = \ABS{z}^2 / (2\tau)$.

    \item[$\ABS{z} \geq \tau$:]
    Now the critical point $\bar\alpha = \ABS{z} - \tau$ is valid, and the positivity of the second derivative $v''(\alpha) = 1/\tau$ classifies this point as a local minimum.
    The corresponding function value is
    \begin{equation*}
      v(\bar \alpha) = \ABS{z} - \tau + \frac{1}{2\tau}\, \tau^2 = \ABS{z} - \frac{\tau}{2}.
    \end{equation*}
  \end{enumerate}

  \noindent
  Finally, for $1 < p < 2$, there is no general closed form solution to \Cref{eq:optim:moreau_implicit}, except for $z = 0$ with solution $\alpha = 0$.
  However, we can use that $v''$ is positive, which implies that any critical point is a local minimum.
  Further, if $z \neq 0$, then it can be readily deduced that the left-hand side of \Cref{eq:optim:moreau_implicit} is smaller than $\ABS{z}$ for $\alpha = 0$ and larger than $\ABS{z}$ for $\alpha = \ABS{z}$.
  Therefore, by continuity it follows that there must be a solution $0 < \bar\alpha < \ABS{z}$, which must be the global minimum due to the global convexity of $v$.
  To simplify the expression for the minimum value, we multiply \Cref{eq:optim:moreau_implicit} with $\alpha$ and solve for $\alpha^p$:
  \begin{equation*}
    \alpha^p = \frac{\alpha(\ABS{z} - \alpha)}{\tau p}.
  \end{equation*}
  Inserting this identity into $v$ with $\alpha = \bar\alpha$ finally yields
  \begin{equation*}
    v(\bar\alpha) = \frac{\alpha(\ABS{z} - \alpha)}{\tau p} + \frac{1}{2\tau}(\ABS{z} - \alpha)^2 = (\ABS{z} - \alpha) \frac{2\alpha + p(\ABS{z} - \alpha)}{2\tau p}.
  \end{equation*}
\end{myproof}
\noindent
The proximal operator can now be easily deduced from the above proof.

\begin{mycorollary}
  The proximal operator $\PROX_{\tau \MODULAR} \colon \SPCX \to \SPCX$ is given as
  \begin{equation}
    \PROX_{\tau \MODULAR}(f)(x) = U_\tau\big( f(x), p(x) \big) \frac{f(x)}{\ABS{f(x)}}
    \label{eq:optim:modular_prox}
  \end{equation}
  with the interpretation $\PROX_{\tau \MODULAR}(f)(x) = 0$ for $f(x) = 0$, where
  \begin{equation}
    U_\tau(z, p) \DEFEQ
    \begin{cases}
      \max\SET{\ABS{z} - \tau,\ 0}, & \text{if } p = 1, \\
      \frac{\ABS{z}}{1 + 2\tau},    & \text{if } p = 2, \\
      \bar\alpha(\ABS{z}, p, \tau)  & \text{else},
    \end{cases}
    \label{eq:optim:modular_prox_mult_func}
  \end{equation}
  and $\bar\alpha$ is defined as in \Cref{lemma:optim:modular_moreau}.
  \label{coroll:optim:modular_prox}
\end{mycorollary}

\begin{myproof}
  The values of the function $U_\tau$ correspond precisely to the minima as calculated in the proof of \Cref{lemma:optim:modular_moreau}.
  For the directional part $f(x) / \ABS{f(x)}$, it suffices to observe that in \Cref{eq:optim:moreau_sphere_min},
  \begin{equation*}
    \inf_{\theta \in S^{m-1}} -\alpha\INNER{z}{\theta}_{\RR^m} = -\alpha\ABS{z}
  \end{equation*}
  is achieved in the unique vector $\bar \theta = z / \ABS{z}$.
\end{myproof}

Often optimization algorithms require the proximal operator of the convex conjugate of a given functional.
Thanks to the Moreau decomposition, this operator can be deduced from the proximal of the functional itself.

\begin{mycorollary}
  The proximal operator $\PROX_{\tau \MODULAR^*}$ of the convex conjugate is given by
  \begin{equation}
    \PROX_{\tau \MODULAR^*}(f)(x) = V_\tau\big( f(x), p(x) \big) \frac{f(x)}{\ABS{f(x)}}
    \label{eq:optim:modular_prox_cconj}
  \end{equation}
  with
  \begin{equation}
    V_\tau(z, p) \DEFEQ
    \begin{cases}
      \min\SET{\ABS{z},\ 1},                            & \text{if } p = 1, \\
      \frac{2\ABS{z}}{\tau + 2},                        & \text{if } p = 2, \\
      \ABS{z} - \bar\alpha(\ABS{z}, p, \tau^{1-p}),     & \text{else}.
    \end{cases}
    \label{eq:optim:modular_prox_cconj_mult_func}
  \end{equation}
  As before, $\bar \alpha = \bar\alpha(\ABS{z}, p, \tau) $ is implicitly defined as solution to \Cref{eq:optim:moreau_implicit}.
  \label{coroll:optim:modular_prox_cconj}
\end{mycorollary}

\begin{myproof}
  Since $\OPS = \OPS^{**}$ for all $\OPS \in \Gamma_0$, we can write the Moreau decomposition \cite[Theorem 14.3]{bauschke_convex_2011} for $\MODULAR^*$ as
  \begin{equation*}
    f = \PROX_{\tau \MODULAR^*}(f) + \tau \PROX_{\tau^{-1} \MODULAR}(\tau^{-1} f).
  \end{equation*}
  This implies that
  \begin{equation*}
    \PROX_{\tau \MODULAR^*}(f) = f - \tau \PROX_{\tau^{-1} \MODULAR}(\tau^{-1} f) = \Big[ \ABS{f} - \tau V_{\tau^{-1}}(\tau^{-1} f(\cdot), p(\cdot) \big) \Big] \frac{f}{\ABS{f}},
  \end{equation*}
  hence by the substitutions $z=f(x)$ and $p=p(x)$, the term in brackets can for $p=1$ be written as
  \begin{equation*}
    \ABS{z} - \tau \max\SET{\tau^{-1} \ABS{z} - \tau^{-1},\ 0} = \ABS{z} - \max \SET{\ABS{z} - 1,\ 0} = \min\SET{\ABS{z},\ 1}.
  \end{equation*}
  Likewise, for $p=2$, we have
  \begin{equation*}
    \ABS{z} - \frac{\ABS{z}}{1 + 2\tau^{-1}} = \ABS{z} \left( 1 - \frac{\tau}{\tau + 2} \right) = \ABS{z}\, \frac{2}{\tau + 2}
  \end{equation*}
  for the bracketed expression.
  In the case $1 < p < 2$, we apply the same substitution to \Cref{eq:optim:moreau_implicit} and get
  \begin{equation*}
    \tau\alpha + p \alpha^{p - 1} = \ABS{z}\quad \Leftrightarrow \quad \beta + \tau^{1 - p} p \beta^{p - 1} = \ABS{z}
  \end{equation*}
  using the new variable $\beta = \tau \alpha$.
  This equation is again of the form \Cref{eq:optim:moreau_implicit} with $\tau^{1 - p}$ in lieu of $\tau$.
\end{myproof}

\paragraph{Solving the implicit equation}
Since the Moreau envelope and the proximal operators involve solving a non-linear equation for each point where $p$ is between 1 and 2, we formulate and analyze in the following a numerical scheme for its solution.
Newton's method seems well-suited for this task since the function whose zero is to be determined is smooth and strictly convex.
The largest practical issue is the determination of a suitable start point that guarantees valid iterates and provides fast convergence to the solution.

\begin{mylemma}
  The Newton iteration for solving \Cref{eq:optim:moreau_implicit} takes the form
  \begin{equation}
    \alpha_{k+1} = \frac{\ABS{z} - p(2-p) \tau \alpha_k^{p-1}}{1 + \tau p(p-1) \alpha_k^{p-2}},\quad \alpha_0 \text{ given},\ k=0,1,\dots.
    \label{eq:optim:modular_newton_iter}
  \end{equation}
  If the start value $\alpha_0$ satisfies
  \begin{equation}
    0 < \alpha_0^{p-1} < \frac{\ABS{z}}{\tau p (2-p)},
    \label{eq:optim:modular_newton_startval_cond}
  \end{equation}
  the sequence $(\alpha_k)_{k \geq 1}$ is positive and strictly increasing, and converges quadratically to the unique root $\bar \alpha$.
  \label{lemma:optim:modular_newton_iter}
\end{mylemma}

\begin{myproof}
  We aim at finding a root of the function $s(\alpha) = \alpha + \tau p \alpha^{p-1} - \ABS{z}$ whose derivative is given by $s'(\alpha) = 1 + \tau p(p-1) \alpha^{p-2}$.
  As noted already in the proof of \Cref{lemma:optim:modular_moreau}, $s$ is strictly increasing and has exactly one root between $\alpha = 0$ and $\alpha = \ABS{z}$, provided that $z \neq 0$.

  For the Newton iterates, we have
  \begin{align*}
    \alpha_{k+1}
    &= \alpha_k - \frac{s(\alpha_k)}{s'(\alpha_k)} = \frac{\alpha_k(1 + \tau p(p-1) \alpha_k^{p-2}) - \alpha_k - \tau p \alpha_k^{p-1} + \ABS{z}}{1 + \tau p(p-1) \alpha_k^{p-2}} \\
    &= \frac{\ABS{z} + \tau p(p-2) \alpha_k^{p-1}}{1 + \tau p(p-1) \alpha_k^{p-2}}.
  \end{align*}
  The first iterate $\alpha_1$ is positive if the numerator $\ABS{z} + \tau p(p-2) \alpha_0^{p-1}$ is positive.
  This is equivalent to \Cref{eq:optim:modular_newton_startval_cond}.

  We now show that the sequence $\alpha_k$ is strictly increasing from $k = 1$ on, which ensures positivity in case $\alpha_1$ is positive.
  The function $s$ is twice differentiable, strictly concave and strictly increasing.
  Using Taylor expansion around the unique root $\bar \alpha$, we have
  \begin{equation*}
    0 = s(\bar \alpha) = s(\alpha_k) + s'(\alpha_k)(\bar \alpha - \alpha_k) + \frac{s''(\eta_k)}{2}(\bar \alpha - \alpha_k)^2
  \end{equation*}
  for some $\eta_k$ between $\bar \alpha$ and $\alpha_k$.
  Dividing by $s'(\alpha_k) > 0$ and rearranging terms yields
  \begin{equation*}
    0 = \underbrace{\frac{s(\alpha_k)}{s'(\alpha_k)} - \alpha_k}_{=-\alpha_{k+1}} + \bar \alpha + \frac{s''(\eta_k)}{2s'(\alpha_k)}(\bar \alpha - \alpha_k)^2,
  \end{equation*}
  i.e. the approximation error $e_k = \bar\alpha - \alpha_k$ satisfies
  \begin{equation*}
    e_{k+1} = -\frac{s''(\eta_k)}{2s'(\alpha_k)}\, e_k^2.
  \end{equation*}
  Besides the quadratic convergence of Newton's method, this also shows that $e_{k+1} > 0$ for all $k \geq 0$ since $s'' < 0$ and $s' > 0$ everywhere.
  Finally, from the iteration we see that
  \begin{equation*}
    \alpha_{k+1} - \alpha_k = -\frac{s(\alpha_k)}{s'(\alpha_k)} > 0
  \end{equation*}
  since $\alpha_k$ is smaller than the root $\alpha$, and $s' > 0$.
  This concludes the proof.
\end{myproof}

\Cref{lemma:optim:modular_newton_iter} formulates a necessary condition that a start value has to meet in order to produce a valid sequence of iterates.
However, the bound does not imply any practical rule how to choose the value and may, in fact, yield impractically large or small numbers.
This situation occurs when $p$ is close to 1 since the resulting power $1 / (p-1)$ applied to the right-hand side of \Cref{eq:optim:modular_newton_startval_cond} is very large and yields either huge or tiny numbers, none of which are useful.

To find a better strategy for the selection of a start value, we repeat the observation from the proof of \Cref{eq:optim:modular_moreau} that $\bar\alpha < \ABS{z}$, which will act as a second bound on the start value.
Finally, since the solutions $\alpha^{(1)}$ and $\alpha^{(2)}$ for the extreme cases $p=1$ and $p=2$, respectively, can be computed explicitly, and since $\alpha^{(2)}$ is positive for $z \neq 0$, we choose the convex combination
\begin{equation}
  \widetilde{\alpha}_0 \DEFEQ (2-p) \alpha^{(1)} + (p-1) \alpha^{(2)}
  \label{eq:optim:modular_newton_startval_guess}
\end{equation}
as an initial guess for a start value and set
\begin{equation}
  \alpha_0 \DEFEQ \min\SET*{\widetilde{\alpha}_0, \ABS{z}, \frac{1}{2}\left(\frac{\ABS{z}}{\tau p (2-p)}\right)^{\frac{1}{p-1}}}
  \label{eq:optim:modular_newton_startval_choice}
\end{equation}
as the effective start value.

\section{TV regularization with variable exponent}
\label{sec:tvreg}

\begin{figure}
  \centering
  \hspace*{3pt}%
  \includegraphics[width=0.47\textwidth]{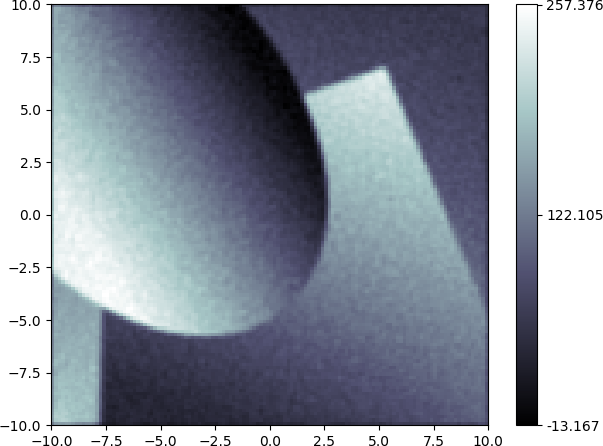}%
  \hspace*{4pt}%
  \includegraphics[width=0.45\textwidth]{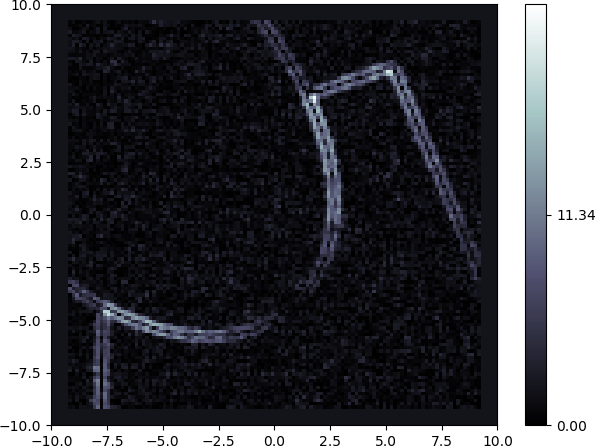}\\[1ex]
  \includegraphics[width=0.45\textwidth]{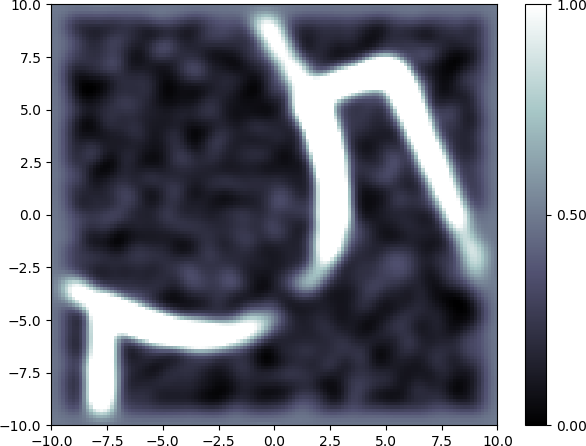}%
  \hspace*{12pt}%
  \includegraphics[width=0.45\textwidth]{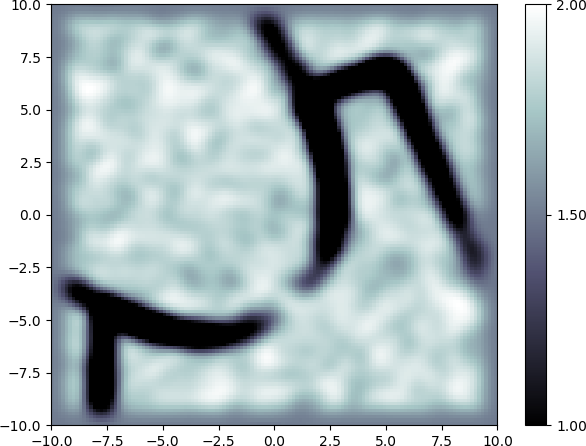}
  \caption{%
  Steps in the generation of the exponent function.
  \emph{Top left:} Input to the exponent computation, a filtered back-projection reconstruction from the low-noise dataset in \Cref{subsec:num:bimodal}.
  \emph{Top right:} Result after taking the absolute value of the smoothed Laplacian of the input.
  \emph{Bottom left:} By smoothing the absolute value and thresholding, the edge areas become broader.
  \emph{Bottom right:} The final exponent function is obtained by subtracting the previous result from 2.
  }
  \label{fig:exp:example}
\end{figure}

Let $\FWDOP \colon \SPCX \to \SPCY$ be a bounded linear operator, where $\SPCX = L^2(\Omega)$ for $\Omega \subset \RR^d$ and $\SPCY$ is another Hilbert space.
We consider the inverse problem of finding $f \in \SPCX$ such that $\FWDOP(f) = g$ for a given $g \in \SPCY$.
TV regularization in the context of $L^2$ spaces can be formulated as the problem of finding
\begin{equation}
  f_\lambda = \ARGMIN_{f \in \SPCX} \left[ \DATAFUNC\big(\FWDOP(f); g\big) + \lambda \TV(f) \right],\quad \lambda > 0,
  \label{eq:tvreg:classical}
\end{equation}
with a given Fr\'{e}chet-differentiable \emph{data discrepancy functional} $\DATAFUNC(\cdot; g)$ and the \emph{total variation}
\begin{equation}
  \TV(f) \DEFEQ
  \begin{cases}
    \NORM{\nabla f}_{L^1(\Omega, \RR^d)},       & \text{if } f \in H^1(\Omega), \\
    +\infty,                                    & \text{else.}
  \end{cases}
\end{equation}
We now replace the 1-norm with the $p$-modular for a given exponent function $p \colon \Omega \to \CINTERV{1}{2}$, leading to the optimization problem
\begin{equation}
  \TVMINL_\lambda \DEFEQ \ARGMIN_{f \in \SPCX} \left[ \DATAFUNC\big(\FWDOP(f); g\big) + \lambda \TV^p(f) \right],
  \label{eq:tvreg:modular}
\end{equation}
where the $p$-TV functional is defined as $\TV^p = \MODULAR \circ \nabla$.
The functional $\TV^p$ imposes two kinds of prior knowledge onto a given inverse problem to determine $f$ from $g = \FWDOP(f)$.
Firstly, it locally enforces either a sparse gradient ($p = 1$) or a smooth solution ($p = 2$), depending on the exponent value at the given location.
Secondly, the exponent itself contains prior knowledge as to which of the two constraints is appropriate to enforce at a specific location.
Hence it is important to find a robust and reliable way to construct an exponent function from given data.
This aspect is investigated further in \Cref{sec:exp}.

\section{Construction of the exponent function}
\label{sec:exp}

\begin{figure}
  \centering%
  \includegraphics[width=0.48\textwidth]{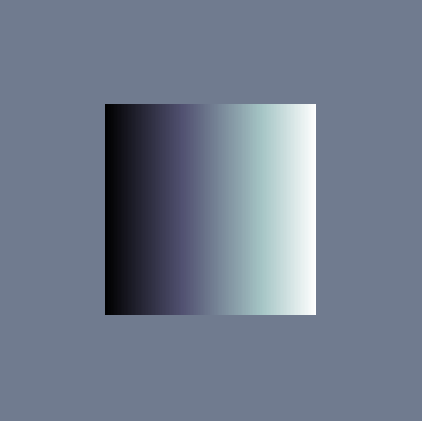}%
  \hspace*{5pt}%
  \includegraphics[width=0.48\textwidth]{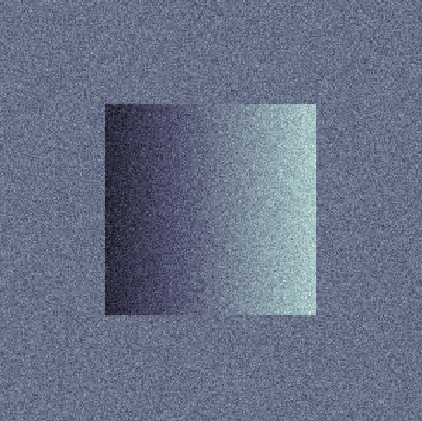}\\
  \vspace*{5pt}%
  \includegraphics[width=0.48\textwidth]{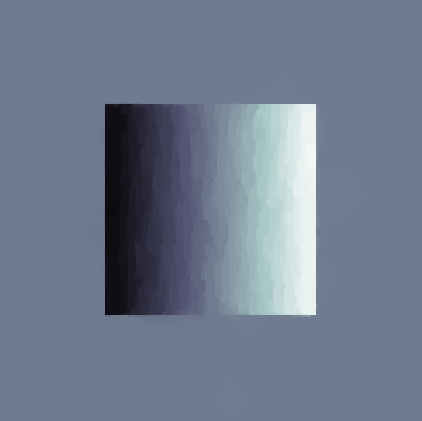}%
  \hspace*{5pt}%
  \includegraphics[width=0.48\textwidth]{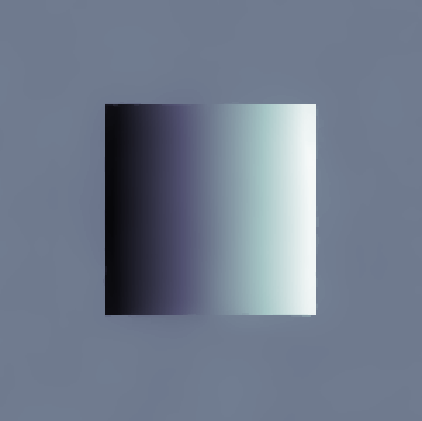}\\
  \vspace*{5pt}%
  \includegraphics[width=0.48\textwidth]{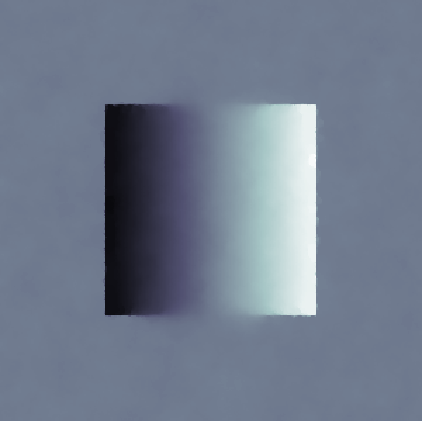}%
  \hspace*{5pt}%
  \includegraphics[width=0.48\textwidth]{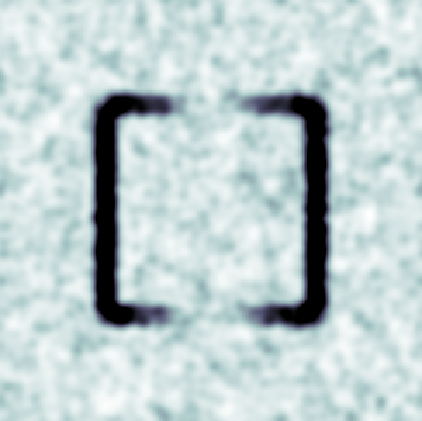}
  \caption{%
  Denoising of a simple image with edges and gradual intensity variation.
  \emph{Top left:} Original image.
  \emph{Top right:} Noisy image.
  \emph{Middle left:} TV denoising result.
  \emph{Middle right:} TGV denoising result.
  \emph{Bottom left:} $\TV^p$ denoising result.
  \emph{Bottom right:} Computed exponent map, where black corresponds to $p = 1$ and white to $p = 2$.
  }
  \label{fig:num:square_denoising}
\end{figure}

In the formulation of the problem \Cref{eq:tvreg:modular}, the exponent $p \colon \Omega \to \CINTERV{1}{2}$ is assumed to be known.
For a concrete application this means that the function needs to be determined beforehand by a method providing a reasonable distinction between edge features corresponding to $p = 1$ and non-edges that are mapped to other exponent values.
Due to numerical instability, it is not feasible to allow exponents arbitrarily close to 1, as seen in \Cref{eq:optim:modular_newton_startval_cond}.

In general, given a function $f \in \SPCX$, we construct the exponent function in the following way:
\begin{enumerate}
  \item
  Compute the smoothed Laplacian $l = \LAPLACE(G_{\sigma_1} \ast f)$, where $G_\sigma$ is a Gaussian of width $\sigma$.
  The smoothing kernel should suppress noise and be chosen rather narrow.
  \item
  Take the absolute value and smooth again to make the region of detected edges larger: $a = G_{\sigma_2} \ast \ABS{l} \colon \Omega \to \RINTERV{0}{\infty}$.
  Typically, this convolution kernel should be significantly wider, i.e., $\sigma_2 > \sigma_1$.
  \item
  Multiply the result with a constant $c$ and threshold at 1, thus effectively clipping at the value $c$: $t = \min\SET{c \cdot a,\ 1} \colon \Omega \to \CINTERV{0}{1}$.
  \item
  The final exponent is given by $p = 2 - c \colon \Omega \to \CINTERV{1}{2}$.
\end{enumerate}

\noindent
We consider two scenarios of generating the exponent prior from given data: ``bootstrapping'' and bimodal imaging.

\paragraph{Bootstrapping}
This technique is applied when only a single dataset is available.
In this case, the first step consists in computing an image from the data that can be used in the procedure for extracting an exponent as shown in \Cref{fig:exp:example}.
For denoising problems, the data is already in image space, and no action is required.
In case of tomography, a simple method like filtered back-projection can be used to gain the initial image.
After that, the procedure as described above is applied to the image, yielding an exponent function that encodes likely occurrences of edges in the prior image.
This step can be enhanced by directly reconstructing edge images by Lambda tomography \cite{faridani_local_1992,quinto_singularities_1993} or Approximate Inverse for feature reconstruction \cite{louis_diffusion_2010,louis_feature_2011}.

The exponent acquired in this way can then be used to define the variable $L^p$ modular that serves as a building block for the $\TV^p$ regularizer in the minimization problem \Cref{eq:tvreg:modular}.

\paragraph{Bimodal imaging}
A regular situation in bimodal imaging, i.e., when the same object is imaged with two different modalities, is that one type of data is much more reliable than the other in terms of signal-to-noise ratio, resolution or similar meaures of quality.

Typical examples are PET-MRI \cite{ehrhardt_joint_2015}, PET-CT \cite{kinahan_dual_2006} and SPECT-CT \cite{bal_combined_2011}.
In all these applications, the actual quantity of interest, a radioactive tracer distribution, cannot be reconstructed without a secondary source of information about the attenuation of the object.
This secondary parameter can be obtained at high resolution, while the primary channel is usually characterized by strong noise and poor spatial resolution.
Nevertheless, edges from the secondary channel may still be valuable prior information for the reconstruction of the primary quantity of interest, but should not be strictly enforced.

The variable exponent in the $\TV^p$ prior seems to be well-suited for this scenario -- it suggests edges without enforcing them.
Here, a reconstruction of the secondary quantity is used as input $f$ to the exponent calculation, and the resulting function $p$ is then used in the regularized problem \Cref{eq:tvreg:modular} for the primary channel alone.

\begin{myremarks}~
  \begin{enumerate}
    \item
    While the exponent prior is not very prone to creating new edges in the primary channel, it may very well lead to edges that are only present in the primary channel being missed by the proposed method.
    Therefore, it should likely be complemented by another reconstruction method that is not biased in this sense.

    \item
    The recipe for generating the exponent $p$ contains 3 parameters $\sigma_1$, $\sigma_2$ and $c$, which may seem like an impractically large parameter space.
    However, there are two aspects which greatly simplify parameter selection, namely
    \begin{inparaenum}[(1)]
      \item that the range of sensible choices is very limited and can be made solely based on the size of the image domain and the level of noise, and
      \item that different settings can be tested with \emph{immediate} and \emph{visual} feedback, due to the simplicity of the involved operations.
    \end{inparaenum}
    This is confirmed by the fact that the examples in \Cref{sec:num} use nearly identical settings.
  \end{enumerate}
\end{myremarks}

\section{Numerical results}
\label{sec:num}

\begin{figure}
  \centering%
  \includegraphics[width=0.48\textwidth]{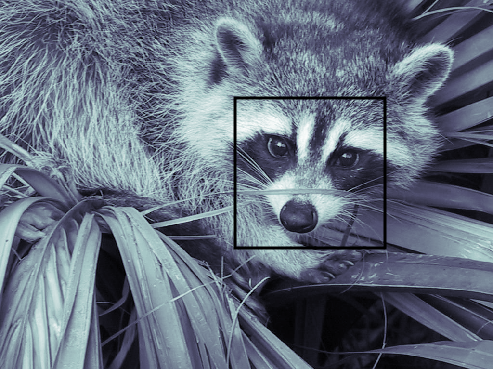}%
  \hspace*{5pt}%
  \includegraphics[width=0.48\textwidth]{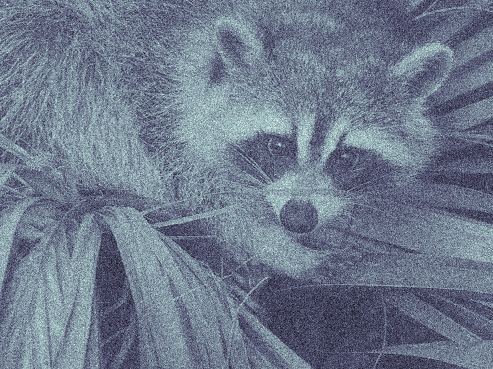} \\
  \vspace*{5pt}%
  \includegraphics[width=0.48\textwidth]{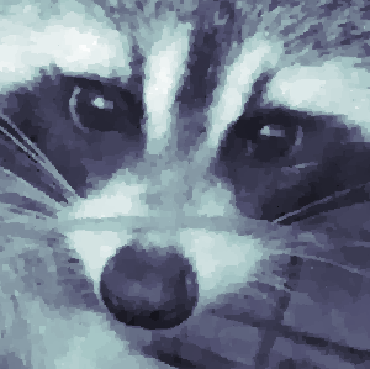}%
  \hspace*{5pt}%
  \includegraphics[width=0.48\textwidth]{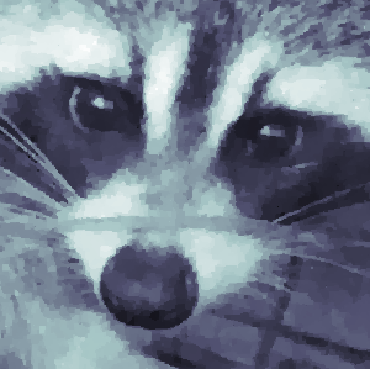}\\
  \vspace*{5pt}%
  \includegraphics[width=0.48\textwidth]{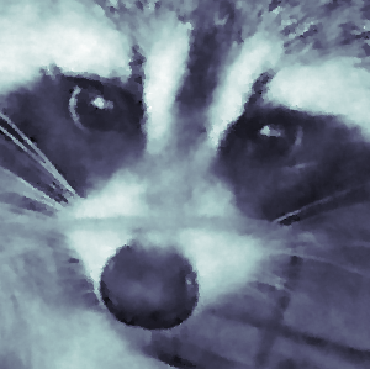}%
  \hspace*{5pt}%
  \includegraphics[width=0.48\textwidth]{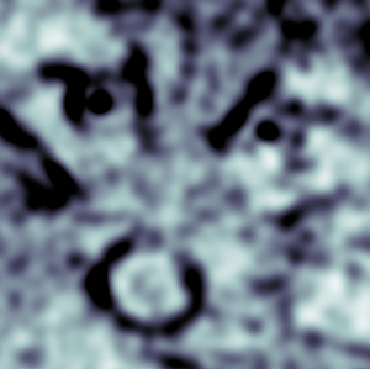}
  \caption{%
  Denoising of a natural image, detail view.
  \emph{Top left:} Original image.
  \emph{Top right:} Noisy image.
  \emph{Middle left:} TV denoising result.
  \emph{Middle right:} TGV denoising result.
  \emph{Bottom left:} $\TV^p$ denoising result.
  \emph{Bottom right:} Computed exponent map, where black corresponds to $p = 1$ and white to $p = 2$.
  }
  \label{fig:num:face_reco}
\end{figure}

In this section we study the numerical characteristics of the proposed regularization functional.
First we demonstrate that structure of the variable $L^p$ modular itself and its derived properties allow a very efficient implementation, thus allowing the functional to be applied in large-scale inverse problems.
We then test the effect of the $\TV^p$ functional in isolation by using it in several single-channel image denoising problems.
Finally, the method of applying the exponent prior from a secondary channel to regularize a primary channel is evaluated in a simulated bimodal tomography example.
The code for all examples is available at GitHub\footnote{\url{https://github.com/kohr-h/odl/tree/variable_lp}}.

\subsection{Implementation}
\label{subsec:num:impl}

The code for the evaluation of the involved functionals and operators is implemented in the ODL (Operator Discretization Library) framework \cite{adler_operator_2017} for inverse problems.
It is written in the Python programming language and features classes and functions for handling vector spaces, operators and functionals in a transparent and efficient way, offering a syntax close to mathematical notation.
The library also contains a number of optimization methods that are adequate for solving convex problems with non-differentiable functionals and compositions with linear operators, among others the Chambolle-Pock method \cite{chambolle_first-order_2010} and two methods based on Forward-Backward splitting \cite{bot_convergence_2015} and Douglas-Rachford splitting \cite{bot_douglas--rachford_2013}, respectively.
Hence, the implementation of the functionals $\rho_p$ and $\rho_p^*$ along with their proximal operators is sufficient to enable numerical tests with regularization using the $\TV^p$ functional.

For running time comparison, the proximal operators are implemented in four variants:
\begin{inparaenum}[(1)]
  \item with \emph{NumPy}\footnote{\url{http://www.numpy.org/}} \cite{van_der_walt_numpy_2011}, a package for fast vectorized array computations,
  \item in the \emph{Cython}\footnote{\url{http://cython.org/}} language \cite{behnel_cython:_2011}, a superset of Python with types that is compiled as native C code,
  \item using \emph{Numba}\footnote{\url{http://numba.pydata.org/}}, a just-in-time compiler to accelerate pure Python functions, and
  \item as GPU kernels using \emph{libgpuarray}\footnote{\url{https://github.com/Theano/libgpuarray}}, a subproject of the deep learning framework Theano \cite{al-rfou_theano:_2016} providing NumPy-like arrays on the GPU.
\end{inparaenum}

\begin{table}[h]
  \begin{tabular}{l|cccccc}
                & NP            & CY            & NB CPU        & NB par        & NB CUDA       & GPUArr        \\
    \hline
    prox        & 1.33          & 0.95 (1.4)    & 1.06 (1.3)    & 0.29 (4.5)    & 0.39 (3.4)    & 0.028 (47)    \\
    cc          & 0.16          & 0.20 (0.8)    & 0.17 (0.9)    & 0.08 (1.8)    & 0.18 (0.86)   & 0.006 (26)
  \end{tabular}
  \caption{%
  Speed comparison of computing the proximal factor $U_\tau$ from \Cref{eq:optim:modular_prox_mult_func} and the convex conjugate integrand $R$ as defined in \Cref{eq:optim:modular_cconj_integrand} for $10^6$ points using 10 Newton iterations.
  The test machine has an 8-core Intel Core i7-6700K CPU and an NVidia GeForce GTX 1070 graphics card.
  The tested implementations are NumPy (NP), Cython (CY), Numba (NB) with CPU, parallel and CUDA targets, and libgpuarray.
  Numbers in parentheses are speed-ups compared to NumPy.
  Software versions: Python 3.5, NumPy 1.12, Cython 0.25.2, Numba 0.31.0, libgpuarray 0.6.0.
  }
  \label{tab:num:speed_comparison}
\end{table}

\noindent
\Cref{tab:num:speed_comparison} shows a comparison of running times of two algorithms that are relevant for regularization with the variable $L^p$ modular.
The first one is the computation of the proximal factor $U_\tau$ as defined in \Cref{eq:optim:modular_prox_mult_func}, which involves a Newton iteration as detailed in \Cref{lemma:optim:modular_newton_iter}.
In the second test, the integrand $R$ of the convex conjugate as given in \Cref{eq:optim:modular_cconj_integrand} is computed for the same number of points.
Clearly, the GPU implementation gives a significant speed-up in both cases ($47 \times$/$26 \times$), while the other acceleration techniques only provide moderate speed gains in the first scenario ($1.4 \times$ to $4.5 \times$) or even result in slower execution in the second case.
This behavior can mainly be attributed to the fact that the function $R$ is too simple in structure and requires too little work for those techniques to pay off, while the evaluation of $U_\tau$ involves a Newton iteration that clearly benefits from acceleration.

\subsection{Denoising examples}
\label{subsec:num:denoising}

\begin{figure}
  \centering
  \includegraphics[width=0.32\textwidth]{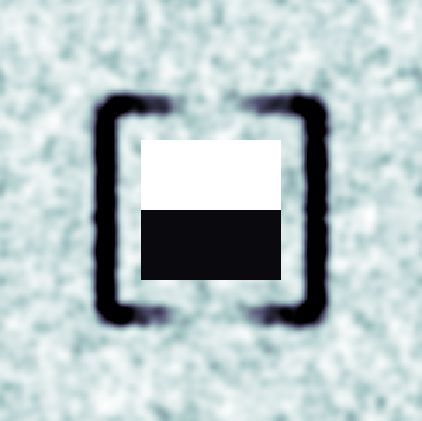}%
  \hspace*{2pt}%
  \includegraphics[width=0.32\textwidth]{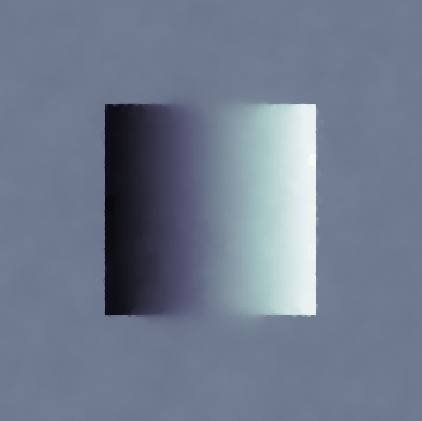}%
  \hspace*{2pt}%
  \includegraphics[width=0.32\textwidth]{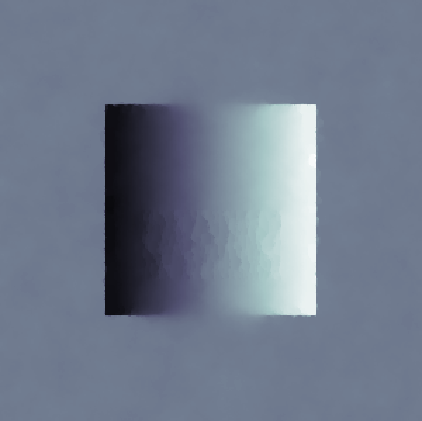}%
  \caption{%
  Effect of variations in the exponent map.
  \emph{Left:} Exponent with new values assigned to selected regions.
  \emph{Middle:} Result when choosing $p = 1.5$ for the lower stripe and $p = 2$ for the upper.
  \emph{Right:} Result when the lower stripe has value $p = 1.05$.
  }
  \label{fig:num:dist_exp}
\end{figure}

We start out by testing $\TV^p$ regularization in a pure denoising context, to keep concerns separate from issues with forward operators.
The first test case is an image whose grey values given as
\begin{equation*}
  f(x) = \chi_{[-5, 5]^2}(x) \cdot x_0
\end{equation*}
on the domain $\Omega = [-10, 10]^2$.
It serves as a simple prototype of an image with sharp edges and gradual intensity variations.
The noisy data is generated by adding standard white noise scaled by 0.1 times the dynamic range of the image.
We will henceforth refer to this as ``10 \%\ noise''.

\begin{figure}
  \centering%
  \begin{minipage}{0.48\textwidth}
    \includegraphics[width=\textwidth]{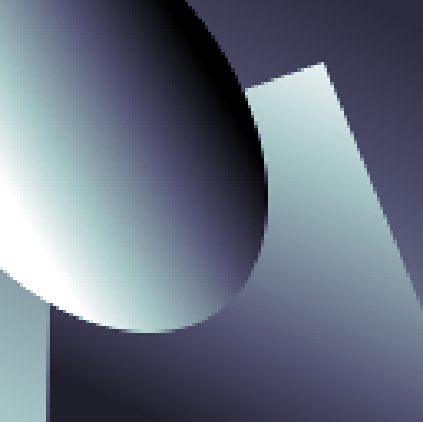}%
  \end{minipage}%
  \hspace*{5pt}%
  \begin{minipage}{0.48\textwidth}
    \includegraphics[width=\textwidth,height=\textwidth]{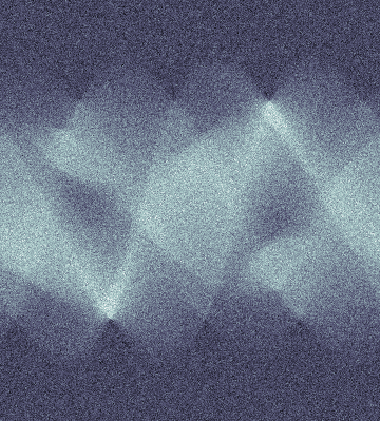}%
  \end{minipage}\\
  \vspace*{5pt}%
  \includegraphics[width=0.48\textwidth]{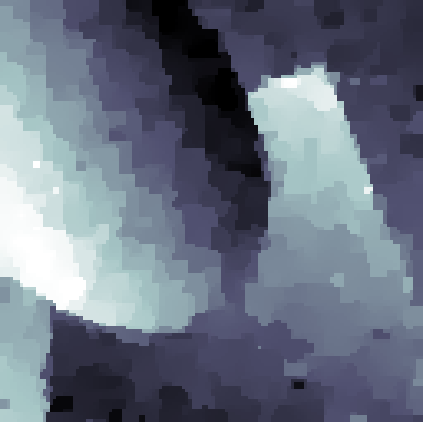}%
  \hspace*{5pt}%
  \includegraphics[width=0.48\textwidth]{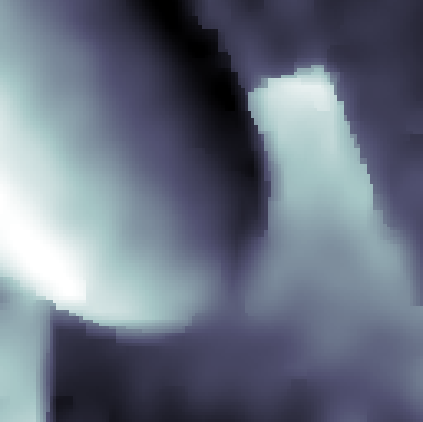}\\
  \vspace*{5pt}%
  \includegraphics[width=0.48\textwidth]{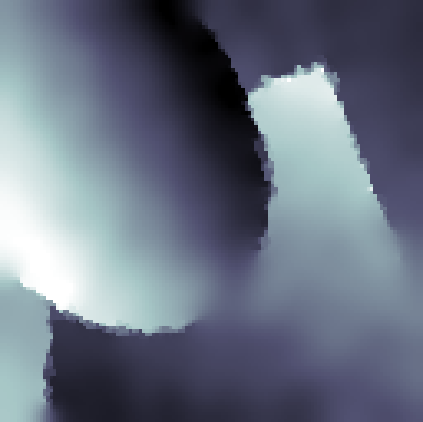}%
  \hspace*{5pt}%
  \includegraphics[width=0.48\textwidth]{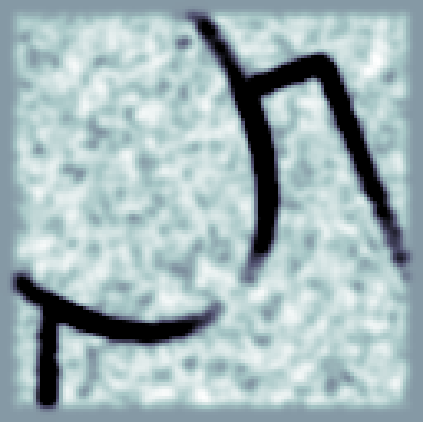}
  \caption{%
  Tomographic reconstruction of an image with edges and gradual intensity variations.
  \emph{Top left:} Original image.
  \emph{Top right:} Noisy data (primary channel).
  \emph{Middle left:} TV reconstruction using only the primary data.
  \emph{Middle right:} TGV reconstruction using only the primary data.
  \emph{Bottom left:} $\TV^p$ reconstruction using the primary data and an exponent map computed from the secondary channel.
  \emph{Bottom right:} Exponent map computed from an FBP reconstruction of the higher-quality secondary data.
  }
  \label{fig:num:bimodal_tomo}
\end{figure}

For comparision, we denoise the image using two other priors, namely classical TV and $\TGV^2$ \cite{bredies_total_2010}.
The latter is implemented by taking the first component $f^*$ of the solution to the product space problem
\begin{equation}
  (f^*, v^*) = \ARGMIN_{(f,v) \in \SPCX \times \SPCX^d} \left[ \NORM{f - g}_2^2 + \lambda_1 \NORM{\nabla f - v}_1 + \lambda_2 \NORM{\mathcal{E}(v)}_1 \right],
  \label{eq:num:tgv}
\end{equation}
with the component-wise gradient $\mathcal{E}(v) = (\nabla v_1, \dots, \nabla v_d) \colon \Omega \to \SPCX^{d\times d}$.

The results in this test scenario are show in \Cref{fig:num:square_denoising}.
Clearly, TGV excels at the given problem since the true image falls into the class of piecewise linear functions that this method tries to promote.
On the other hand, TV exhibits the typical staircasing effect, which makes it a poor fit for this particular phantom.

Looking at the results with $\TV^p$ prior, one can make some interesting observations.
First, the staircasing problem in the interior of the square does not show here since the exponent map takes values between 1.5 and 2.0, resulting in the smoothing behavior that is also typical for the Tikhonov-type prior $\NORM{\nabla f}_2^2$.
Second, moderate variations of the exponent due to noise have little influence on the values of the reconstruction (see also \Cref{fig:num:dist_exp}).
Thirdly, the transition from $p = 1$ to $p > 1$ does not introduce artificial edges.
Any artificial boundaries stem from the TV prior itself that is effective in regions where $p = 1$.

On the downside, the variable $L^p$ regularization blurs the top and bottom edges in the middle where the contrast to the background is low, which is to be expected given the exponent map not capturing these edges.
A more fine-tuned exponent calculation would help alleviate this issue.
Further, the left and right edges appear slightly jagged, compared to the other two alternatives.
This issue can likely be tackled by choosing a slightly larger regularization parameter.

The second test scenario consists in denoising of a natural image with more complex edge structure and intensity variation.
The test image is taken to be the transformation to greyscale of the image generated by the SciPy\footnote{\url{https://scipy.org/}} command \texttt{scipy.misc.face()}, and its noisy variant contains 15 \% white noise.

Results for this test case are shown in \Cref{fig:num:face_reco}.
It is remarkable that the difference between TV and TGV are only marginal.
While the exact reason for this behavior remains unclear, it can be said that the parameter $\lambda_2$ in \Cref{eq:num:tgv} has practically no influence even if varied by 10 orders of magnitude.
Hence, any lowering of the parameter value for $\lambda_1$ -- to reduce staircasing -- below the value used in the TV method results in insufficient noise suppression.

By contrast, the $\TV^p$ approach works very well in this scenario and produces a visually pleasing and naturally looking result with smooth intensity variations in the regions where they ought to be expected, while preserving most of the sharp edges.
The exact behavior of the method can easily be interpreted from the generated exponent map.

\subsection{Bimodal tomography examples}
\label{subsec:num:bimodal}

This application scenario considers tomography with two datasets generated from the same phantom (taken from \cite{bredies_total_2010}).
It is intended to indicate the potential of variable $L^p$ TV regularization in bimodal imaging, not to model realistic imaging conditions.

\begin{figure}
  \centering%
  \includegraphics[width=0.48\textwidth]{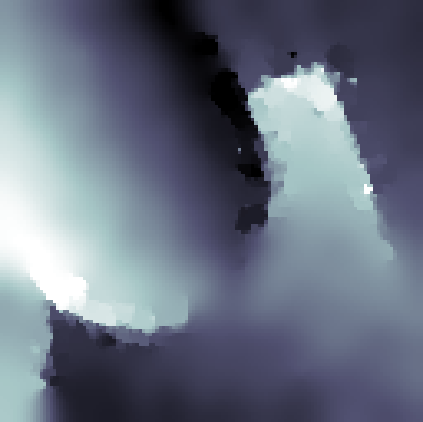}%
  \hspace*{5pt}%
  \includegraphics[width=0.48\textwidth]{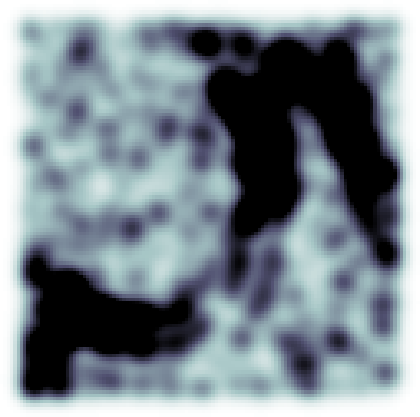}
  \caption{%
  Reconstruction from the same data as in \Cref{fig:num:bimodal_tomo}, but without using the secondary channel.
  \emph{Left:} $\TV^p$ reconstruction.
  \emph{Right:} Exponent map computed from the very noisy dataset using an over-regularized FBP.
  }
  \label{fig:num:bimodal_tomo_varlp_bootstrap}
\end{figure}

The datasets are the 2D divergent beam transform of a simple digital phantom corrupted by low noise (1 \%) in the secondary and high noise (15 \%) in the primary imaging channel, respectively.
In \Cref{fig:num:bimodal_tomo}, both TV and TGV reconstructions make use of the more noisy primary channel only, while $\TV^p$ takes the secondary, less noisy data into account for computing the exponent map.
For a fair comparison, \Cref{fig:num:bimodal_tomo_varlp_bootstrap} shows the result of $\TV^p$ regularization with the ``bootstrapping'' approach of using only the primary channel.
In both cases, the exponent function $p$ is computed from a filtered back-projection reconstruction with the approach described in \Cref{sec:exp}.

As expected, the TV reconstruction recovers the edges well at the cost of strong staircasing within the gradually varying regions.
TGV, on the other hand, removes the staircasing, at the cost of turning some edges into linear slopes.
Clearly, $\TV^p$ yields the visually best result when the secondary channel information is taken into account in form of an accurate exponent map.
Only at the boundaries some parts of the edges are blurred due to the regions of $p = 1$ not extending all the way to the boundaries, a consequence of the convolution used for smoothing.
When taking only the primary channel into account for both exponent and reconstruction, the result contains a larger amount of the blocky structures of the TV reconstruction since the regions with $p = 1$ are overestimated.
More sophisticated methods for finding edges would potentially help alleviate this issue.

\section{Future work}
\label{sec:future}

While the $\TV^p$ prior using the modular $\rho_p$ leads to very encouraging results in the test scenarios considered in this paper, it has the potential drawback of not being scale-invariant.
In other words, the absolute magnitude of the unknown image may not only have an impact on the choice of regularization parameter, but also on the relative importance of regions with different exponents to the functional value.
This issue can be solved by replacing the variable $L^p$ modular with its corresponding norm \Cref{eq:varlp:def_norm}.
Furthermore, the $p$-norm seems to be better suited for the derivation of guarantees for convergence and stability.
This topic, along with efficient ways of computing the norm, its convex conjugate and its proximal operator, will be the subject of future work.

\section*{Acknowledgments}
Part of this work was financially supported by the Netherlands Organization for Scientific Research (NWO), project 639.073.506.

\end{document}